\documentclass[10pt]{article}
\usepackage{latexsym}
\usepackage{amssymb}
\usepackage{epsfig}

\newtheorem{theorem}{Theorem}[section]

\newtheorem{proposition}[theorem]{Proposition}
\newtheorem{corollary}[theorem]{Corollary}
\newtheorem{conjecture}[theorem]{Conjecture}

\newcommand{\eps}{\varepsilon}
\newcommand{\beq}{\begin{equation}}
\newcommand{\eeq}{\end{equation}}
\newcommand{\prob}{{\mathbb P}}
\newcommand{\expec}{{\mathbb E}}

\newcommand{\sss}{\scriptscriptstyle }
\newcommand{\tinbr}[1]{{\scriptscriptstyle (#1)}}

\newcommand{\proof} {\noindent {\bf Proof}. \hspace{2mm}}
\newcommand{\qed}    {\hspace*{\fill} $\Box$ \medskip}

\title{Large deviations for eigenvalues of sample covariance matrices, with applications to mobile communication systems}
\author{Anne Fey
\thanks{Delft Institute of Applied Mathematics, 
Faculty of EEMCS, 
Delft University of Technology,
The Netherlands,
{\tt a.c.fey-denboer@tudelft.nl}}
\and
Remco van der Hofstad$^{\dagger}$
\thanks{
Department of Mathematics,
Eindhoven University,
of Technology,
P.O.\ Box 513,
5600 MB Eindhoven,
The Netherlands,
{\tt rhofstad@win.tue.nl}}
\and Marten J.\ Klok
\thanks{
ORTEC BV,
Orlyplein 145c,
1043 DV Amsterdam,
The Netherlands,
{\tt MKlok@ortec.nl}}}

\begin{document}

\maketitle

\begin{abstract}
We study sample covariance matrices of the form $W=\frac 1n C C^T$, where
$C$ is a $k\times n$ matrix with i.i.d.\ mean zero entries. This is a
generalization of so-called Wishart matrices, where the entries of $C$ are independent
and identically distributed standard normal random variables.
Such matrices arise in statistics as sample covariance matrices,
and the high-dimensional case, when $k$ is large,
arises in the analysis of DNA experiments.

We investigate the large deviation properties of the largest and
smallest eigenvalues of $W$ when either $k$ is fixed and $n\rightarrow \infty$,
or $k_n\rightarrow \infty$ with $k_n=o(n/\log\log{n})$,
in the case where the squares of the i.i.d.\ entries have finite
exponential moments. Previous results, proving a.s.\ limits
of the eigenvalues, only require finite fourth moments.

Our most explicit results for $k$ large
are for the case where the entries of $C$ are $\pm1$
with equal probability. We relate the large deviation rate
functions of the smallest and largest eigenvalue to the rate functions
for independent and identically distributed standard normal entries
of $C$.
This case is of particular interest, since it is related to the
problem of the decoding of a signal in a code division multiple
access system arising in
mobile communication systems. In this example, $k$ plays the role of the number of
users in the system, and $n$ is the length of the coding sequence of
each of the users. Each user transmits at the same time and uses the
same frequency, and the codes are used to distinguish the signals
of the separate users. The results imply large deviation bounds for
the probability of a bit error due to the interference of the various
users.
\end{abstract}


\pagestyle{myheadings} \markboth{\sf Fey-den Boer, van der Hofstad and Klok}
{\sf Large deviations for eigenvalues of sample covariance matrices}

\section{Introduction}

The sample covariance matrix $W$ of a matrix $C$ with $k$ rows and $n$ columns is defined as
$\frac 1n C C^T$. If $C$ has random entries, then the spectrum of $W$ is random as well.
Typically, $W$ is studied in the case that $C$ has i.i.d.\ entries, with mean 0
and variance 1. For this kind of $C$, it is known that when $k,n \to \infty$ such
that $k/n = \beta$, where $\beta$ is a constant,
the eigenvalue density tends to a deterministic density \cite{march}.
The boundaries of the support of this distribution are $(1-\sqrt\beta)_+^2$ and
$(1+\sqrt\beta)^2$, where $x_+=\max\{0,x\}$. This suggests that the smallest
eigenvalue $\lambda_{\rm min}$ converges to $(1-\sqrt\beta)_+^2$, while
the largest eigenvalue $\lambda_{\rm max}$ converges to $(1+\sqrt\beta)^2$.
Bai and Yin \cite{bai1} have proved a.s.\ convergence of $\lambda_{\rm min}$ to $(1-\sqrt\beta)_+^2$.
Bai, Silverstein and Yin \cite{bai2} proved a.s.\ convergence of
$\lambda_{\rm max}$ to $(1+\sqrt\beta)^2$, see also
\cite{YinBaiKri88}. The strongest results
apply in the case that all entries of $C$ are i.i.d.\ with
mean 0, variance 1 and finite fourth moment.
Related results, including a central limit theorem for the
linear spectral statistics, can be found in \cite{BaiSil98, BaiSil04},
to which we also refer for an overview of the extensive literature.

In the special case that the entries of $C$ have a standard normal
distribution, $W$ is called a {\it Wishart matrix}. Wishart matrices
play an important role in multivariate statistics as they describe
the correlation structure in i.i.d.\ Gaussian multivariate data.
For Wishart matrices, the large deviation rate function for
the eigenvalue density with rate $\frac{1}{n^2}$ has been
derived by Guionnet \cite{guion} and Hiai and Petz \cite{hiai}.
However, the proofs depend heavily on the fact that $C$
has standard normal i.i.d.\ entries, for which the density of the
ordered eigenvalues can be explicitly computed.

In this article, we investigate the large deviation
rate functions with rate $\frac 1n$ of the smallest and largest
eigenvalue of $W$, for certain non-Gaussian entries of $C$.
We pose a strong condition on the tails of the
entries, by requiring that the exponential moment
of the square of the entries is bounded in a
a neighborhood of the origin. We shall also comment
on this assumption, which we believe to be necessary
for our results to apply.

We let $n \to \infty$, and $k$ is either fixed or tends to infinity
not faster than $o(n/\log\log n)$.
Our results imply that all eigenvalues tend to 1 and that all other
values are large deviations. We obtain the asymptotic
large deviation rate function of $\lambda_{\rm min}$ and
$\lambda_{\rm max}$ when $k\rightarrow \infty$.
In certain special cases, we can show that
the asymptotic large deviation rate function is equal to the one
for Wishart matrices, which can be interpreted as saying that
the spectrum of sample covariance matrices with $k$ and $n$ large
is close to the one for i.i.d.\ standard normal entries.
This proves a kind of universality result for the large deviation
rate functions.

This paper is organized as follows. In Section \ref{general},
we derive an explicit expression for the large deviation rate
functions of $\lambda_{\rm min}$ and $\lambda_{\rm max}$. In Section
\ref{special}, we calculate lower bounds for the case that the
distribution of $C_{mi}$ is symmetric around 0, and $|C_{mi}| < M$
almost surely, for some $M>0$.
In Section \ref{sec-pm1}, we specialize to the case where $C_{mi}=\pm 1$
with equal probability, which arises in an application in wireless
communication. We describe the implications of our results in this application
in Section \ref{cdma}. Part of the results for this application have
been presented at an electrical engineering conference
\cite{fey}.

\section{General mean zero entries of $C$}
\label{general}
In this section, we prove large deviation results for the smallest
and largest eigenvalues of sample covariance matrices.

\subsection{Large deviations for $\lambda_{\rm min}$ and $\lambda_{\rm max}$}
Define $W=\frac 1n C C^T$ to be the matrix of sample covariances.
We denote by $\mathbb{P}$ the law of $C$ and by $\mathbb{E}$ the corresponding
expectation. Throughout the paper, we assume that the i.i.d.\ real matrix elements
of $C$ are normalized, i.e.,
    \beq
    \label{Cass}
    \expec[C_{ij}]=0, \qquad {\rm Var}(C_{ij})=1.
    \eeq
The former implies that a.s., the off diagonal elements of the matrix $W$ converge to
zero, the second implies that the diagonal elements converge to 1, a.s. By a rescaling argument,
the second assumption is without loss of generality.

In this section, we rewrite the probability for a large deviation
of the largest and smallest eigenvalues of $W$, $\lambda_{\rm max}$
and $\lambda_{\rm min}$, respectively, into that of a large deviation
of a sum of i.i.d.\ random variables. This rewrite allows us to use
Cram\'{e}r's Theorem to obtain an expression for the
rate function. This section gives a heuristic derivation of
our result, that will be turned into a proof in Section \ref{a,b}.

For any matrix $W$, and any vector $\mathbf{x}$ with $k$ coordinates
and norm $\|\mathbf{x}\|_2 = 1$, we have
    $$
    \lambda_{\rm min} \leq \langle \mathbf{x}, W\mathbf{x}\rangle \leq  \lambda_{\rm max}.
    $$
Moreover, for the normalized eigenvector $\mathbf{x}_{\rm min}$ corresponding
to $\lambda_{\rm min}$, the lower bound is attained,
while for the normalized $\mathbf{x}_{\rm max}$, corresponding
to $\lambda_{\rm max}$, the upper bound is attained. Therefore, we can write
    \begin{eqnarray}
     P_{\rm min}(\alpha)&
     = \mathbb{P}(\lambda_{\rm min} \leq \alpha)& = \prob(\exists \mathbf{x}: \|\mathbf{x}\|_2 = 1,
    \langle \mathbf{x}, W\mathbf{x}\rangle \leq \alpha),\\
     \nonumber P_{\rm max}(\alpha)
     & = \mathbb{P}(\lambda_{\rm max} \geq \alpha)& = \prob(\exists \mathbf{x}: \|\mathbf{x}\|_2 = 1,
    \langle \mathbf{x}, W\mathbf{x}\rangle \geq \alpha).
    \label{Pmin}
    \end{eqnarray}

We use that the above is the probability of a union of events, and
bound this probability from below by considering only one $\mathbf{x}$,
and from above by summing over all $\mathbf{x}$. Since there are uncountably many
possible $\mathbf{x}$, we will do this approximately by summing over a finite
number of vectors. The lower bound for the probability yields
an upper bound for the rate function, and vice versa.

We first heuristically explain the form of the rate function
of $\lambda_{\rm max}$ and $\lambda_{\rm min}$, and highlight the proof.
The special form of a sample covariance matrix allows us to rewrite
    \beq
    \label{Sxrew}
    \langle \mathbf{x}, W\mathbf{x}\rangle = \frac 1n \|C^T \mathbf{x}\|_2^2
    = \frac 1n \sum_{i=1}^n \left(\sum_{m=1}^k x_m C_{mi}\right)^2
    =\frac 1n \sum_{i=1}^n S_{\mathbf{x},i}^2,
    \eeq
where
    \beq
    S_{\mathbf{x},i} = \sum_{m=1}^k x_m C_{mi},
    \label{Sx}
    \eeq
    with $S_{\mathbf{x},i}$ i.i.d.\ for $ i=1, \ldots, m$.
Define
    \beq
    \label{Ikdef}
    I_k(\alpha) = \inf_{\mathbf{x}\in \mathbb{R}^k:\|\mathbf{x}\|_2=1} \sup_t \left(t\alpha
    -\log \expec[e^{t S_{\mathbf{x},1}^2}]\right).
    \eeq
Since $\mathbb{E}[S_{\mathbf{x},1}^2]=1$, and $t\mapsto \log \expec[e^{t S_{\mathbf{x},1}^2}]$
is increasing and convex, we see that, for fixed $\mathbf{x}$, the optimal $t$ is
non-negative for $\alpha\geq 1$ and non-positive for $\alpha\leq 1$. The sign of $t$ will
play an important role in the proofs in Sections \ref{special}--\ref{sec-pm1}.

We can now state the first result of this paper.

    \begin{theorem}\label{prop-LDeig}
    Assume that (\ref{Cass}) holds. Then,\\
    (a) for all $\alpha\geq 1$ and fixed $k\geq 2$
    \beq
    \label{ratemaxa}
    \limsup_{n\rightarrow \infty}-\frac 1n \log \mathbb{P}(\lambda_{\rm max}
    \geq \alpha) \leq I_k(\alpha),
    \eeq
    and
    \beq
    \label{ratemaxb}
    \liminf_{n\rightarrow \infty}-\frac 1n \log \mathbb{P}(\lambda_{\rm max}
    \geq \alpha) \geq \lim_{\eps\downarrow 0} I_k(\alpha-\eps),
    \eeq
    (b) for all $0 \leq \alpha \leq 1$ and fixed $k\geq 2$
    \beq
    \label{ratemina}
    \limsup_{n\rightarrow \infty}
    -\frac 1n \log \mathbb{P}(\lambda_{\rm min} \leq \alpha) \leq I_k(\alpha),
    \eeq
    and
    \beq
    \label{rateminb}
    \liminf_{n\rightarrow \infty}
    -\frac 1n \log \mathbb{P}(\lambda_{\rm min} \leq \alpha)
    \geq \lim_{\eps\downarrow 0} I_k(\alpha+\eps).
    \eeq
    When there exists an $\epsilon>0$ such that $\expec[e^{\epsilon C_{11}^2}]<\infty$ and when ${\rm Var}(C_{11}^2)>0$,
    then $I_k(\alpha)>0$ for all $\alpha\neq 1$.
    \end{theorem}
We will now discuss the main result in Theorem \ref{prop-LDeig}.
Theorem \ref{prop-LDeig} is only useful when
$I_k(\alpha)>0$, which we prove under the strong
condition that there exists an $\epsilon>0$ such that $\expec[e^{\epsilon C_{11}^2}]<\infty$. For example, a.s.\ limits for the largest and smallest
eigenvalues are proved under the {\it much} weaker
condition that the fourth moment of the matrix entries
$C_{im}$ is finite. However, it is well known that the
exponential bounds present in large deviations are only valid when the
random variables under consideration have finite exponential moments
(see e.g., Theorem \ref{thm-CT} below). In this case,
the rate functions can be equal to zero, and the large deviation
results are rather uninformative. Since the eigenvalues are
{\it quadratic} in the entries $\{C_{im}\}_{i,m}$, this translates
into the above condition, which we therefore believe to be {\it necessary}.

Secondly, we note that,
due to the occurrence of an infimum over $\mathbf{x}$ and a supremum over $t$,
it is unclear whether the function $\alpha\mapsto I_k(\alpha)$ is continuous.
Clearly, when $\alpha\mapsto I_k(\alpha)$ is continuous, the upper and lower bounds in
(\ref{ratemaxa}) and (\ref{ratemaxb}), as well as the ones in
(\ref{ratemina}) and (\ref{rateminb}), are equal. We will see that
this is the case for Wishart matrices in Section \ref{wishart}. The function
$\alpha \mapsto I_k(\alpha)$ can easily be seen to be increasing on
$[1,\infty)$ and decreasing on $(0,1]$, since $\alpha\mapsto
\sup_t \left( t\alpha-\log \expec[e^{t S_{\mathbf{x},1}^2}]\right)$ has
the same monotonicity properties for every fixed $\mathbf{x}$, so that the limits
$\lim_{\eps\downarrow 0} I_k(\alpha+\eps)$ and
$\lim_{\eps\downarrow 0} I_k(\alpha-\eps)$ exist as monotone limits.
The continuity of $\alpha\mapsto I_k(\alpha)$ is not obvious.
For example, in the simplest case where
$C_{ij}=\pm 1$ with equal probability, we know that the large deviation rate function
is {\it not} continuous, since the largest eigenvalue is at most $k$. Therefore,
$\mathbb{P}(\lambda_{\rm max} \geq \alpha)=0$ for any $\alpha>k$,
and, if $\alpha\mapsto I_k(\alpha)$ is the rate function of
$\lambda_{\rm max}$ for $\alpha\geq 1$, then $I_k(\alpha)=\infty$
for $\alpha>k$. It remains an interesting problem to determine in
what cases $\alpha\mapsto I_k(\alpha)$ {\it is} continuous.

Finally, we only prove that $I_k(\alpha)>0$ for all $\alpha\neq 1$ when
${\rm Var}(C_{11}^2)>0$. By the normalization that $\expec[C_{11}]=0,
\expec[C_{11}^2]=1$, this only excludes the case where $C_{11}=\pm 1$
with equal probability. This case will be investigated in
more detail in Theorem \ref{Iuitrekenen3}, where we shall also prove a lower bound
implying that $I_k(\alpha)>0$ for all $\alpha\neq 1$.

Denote
    \beq
    \label{LDext}
    I_k(\alpha,\beta)=\inf_{\stackrel{\mathbf{x, y}\in \mathbb{R}^k:\|\mathbf{x}\|_2=\|\mathbf{y}\|_2=1}
    {\langle \mathbf{x}, \mathbf{y}\rangle=0}}
    \sup_{s,t} \left( t\alpha+s\beta
    -\log \expec[e^{t S_{\mathbf{x},1}^2+s S_{\mathbf{y},1}^2}]\right).
    \eeq
Our proof also reveals that, and for all $0\leq \beta\leq 1, \alpha\geq 1$,
    \beq
    \label{extensiontwoa}
    \limsup_{n\rightarrow \infty}-\frac 1n \log \mathbb{P}(\lambda_{\rm max}
    \geq \alpha, \lambda_{\rm min} \leq \beta) \geq I_k(\alpha, \beta),
    \eeq
and
    \beq
    \label{extensiontwob}
    \lim_{n\rightarrow \infty}-\frac 1n \log \mathbb{P}(\lambda_{\rm max}
    \geq \alpha, \lambda_{\rm min} \leq \beta) \leq \lim_{\eps\downarrow 0}
    I_k(\alpha+\eps, \beta-\eps).
    \eeq
For Wishart matrices, for which the entries of $C$ are i.i.d.\ standard normal,
the random variable $S_{\mathbf{x},i}$ has a standard normal
distribution, so that we can explicitly calculate $I_k(\alpha)$.
We will elaborate on this in Section \ref{wishart} below.
For the case that $C_{mi} = \pm 1$ with equal probabilities, Theorem \ref{prop-LDeig}
and its proof have also appeared in \cite{fey}.

\subsection{Proof of Theorem \ref{prop-LDeig}(a) and (b)}
\label{a,b}
In the proof, we will repeatedly make use of the {\it largest-exponent-wins principle}.
We first give a short explanation of this principle.
This principle is about the exponential rate of the sum of two (or more)
probabilities. From this point, we will abbreviate `exponential rate of a probability'
by `rate'. Because of the minus sign, a smaller rate $I$ means a larger exponent,
and thus a larger probability. Thus, if for two events $E_1$ and $E_2$, both depending on
some parameter $n$, we have
    $$
    \mathbb{P}(E_1)\sim e^{-nI_1} \quad \mbox{ and } \quad \mathbb{P}(E_2)\sim e^{-nI_2}
    $$
then
    \beq
    -\lim_{n\to\infty}
    \frac1n\log(\mathbb{P}(E_1) + \mathbb{P}(E_2)) \sim \min\{I_1, I_2\}.
    \label{largestexponentwins}
    \eeq
In words, the principle states that as $n \to \infty$, the
smallest exponent (i.e., the largest rate) will become negligible.
It also implies that
    \beq
    -\lim_{n\to\infty}
    \frac1n\log\mathbb{P}(E_1\cup E_2) \sim \min\{I_1, I_2\}.
    \label{LEWevents}
    \eeq

In the proof, we will make essential use of Cram\'er's Theorem, which
we state here for the sake of completeness:

\begin{theorem}[Cram\'er's theorem and Chernoff bound]
\label{thm-CT}
Let $(X_i)_{i=1}^{\infty}$ be a sequence of i.i.d.\ random variables. 
Then, for all $a\geq\mathbb{E}[X_1]$,
    \beq
    \label{CTgeq}
    -\lim_{n\to\infty}
    \frac1n\log\mathbb{P}(\frac 1n \sum_{i=1}^n X_i \geq a)
    =\sup_{t\geq 0} \left( ta -\log\mathbb{E}[e^{tX_1}]\right),
    \eeq
while, for all $a\leq \mathbb{E}[X_1]$,
    \beq
    \label{CTleq}
    -\lim_{n\to\infty}
    \frac1n\log\mathbb{P}(\frac 1n \sum_{i=1}^n X_i \leq a)
    =\sup_{t\leq 0} \left(ta -\log\mathbb{E}[e^{tX_1}]\right).
    \eeq
The upper bounds in (\ref{CTgeq})-(\ref{CTleq}) hold for every $n$.\\
Furthermore, when $\mathbb{E}[e^{tX_1}]<\infty$ for all
$t$ with $|t|\leq \epsilon$ and some
$\epsilon>0$, then the right-hand sides of (\ref{CTgeq}) and (\ref{CTleq})
are strictly positive for all $a\neq \mathbb{E}[X_1]$.
\end{theorem}
See e.g.,
\cite[Theorem 1.1, pages 5-6 and Proposition 1.9, page 13]{OliVar05}
for this result, and see \cite{dembo} and \cite{denhollander} for general
introductions to large deviation theory.

For the proof, we start by showing that $I_k(\alpha)>0$
for all $\alpha\neq 1$ when there exists an $\epsilon>0$ such that $\expec[e^{\epsilon C_{11}^2}]<\infty$
and when ${\rm Var}(C_{11}^2)>0$.
For this, we note that, by the Cauchy-Schwarz
inequality and (\ref{Sx}), for every $\mathbf{x}$ with
$\|\mathbf{x}\|_2=1$,
    \[
    S_{\mathbf{x},i}^2\leq \sum_{m=1}^k x_m^2 \sum_{m=1}^k C_{mi}^2=\sum_{m=1}^k C_{mi}^2,
    \]
so that $\expec [e^{t S_{\mathbf{x},i}^2}]
\leq \expec [e^{t C_{11}^2}]^k<\infty$ whenever there exists an $\epsilon>0$ such that $\expec[e^{\epsilon C_{11}^2}]<\infty$.
Thus, uniformly in $\mathbf{x}$ with
$\|\mathbf{x}\|_2=1$, the random variables $S_{\mathbf{x},i}^2$
have bounded exponential moments for $t\leq \epsilon$.
As a result, the Taylor expansion
    \begin{equation}
    \label{TaylorexpMGF}
    \log\expec [e^{t S_{\mathbf{x},i}^2}]
    =t+\frac{t^2}{2}{\rm Var}(S_{\mathbf{x},i}^2)+ {\cal O}(|t|^3)
    \end{equation}
holds uniformly in $\mathbf{x}$ with
$\|\mathbf{x}\|_2=1$.
We compute, since $\mathbb{E}[S_{\mathbf{x},i}^2]=\mathbb{E}[C_{11}^2]=1$,
and for $\mathbf{x}$ with
$\|\mathbf{x}\|_2=1$,
    \[
    \expec[S_{\mathbf{x},i}^4]=3\Big(\sum_{m} x_m^2\Big)^2-3\sum_m x_m^4
    +\mathbb{E}[C_{11}^4]\sum_m x_m^4=3-3\sum_m x_m^4
    +\mathbb{E}[C_{11}^4]\sum_m x_m^4,\]
that
    \[{\rm Var}(S_{\mathbf{x},i}^2)
    =3-3\sum_m x_m^4
    +\mathbb{E}[C_{11}^4]\sum_m x_m^4
    -1=2-2\sum_m x_m^4
    +{\rm Var}(C_{11}^2)\sum_m x_m^4,
    \]
which is bounded, since
by assumption $\expec [e^{t C_{11}^2}]<\infty$.
Furthermore, $\sum_m x_m^4\in [0,1]$ uniformly
in $\mathbf{x}$ with $\|\mathbf{x}\|_2=1$, so that,
again uniformly
in $\mathbf{x}$ with $\|\mathbf{x}\|_2=1$,
${\rm Var}(S_{\mathbf{x},i}^2)\geq
\min\{2, {\rm Var}(C_{11}^2)\}>0$.
We conclude that, for $t$ sufficiently small,
uniformly in $\mathbf{x}$ with
$\|\mathbf{x}\|_2=1$, and by ignoring higher-order
Taylor expansion terms of $t\mapsto \log\expec [e^{t S_{\mathbf{x},i}^2}]$
in (\ref{TaylorexpMGF}), which is allowed when $|t|$ is sufficiently small,
    \[
    \log\expec [e^{t S_{\mathbf{x},i}^2}]\leq
    t+t^2\min\{2, {\rm Var}(C_{11}^2)\}.
    \]
In turn, this implies that for $|t|\leq \epsilon$
small, and uniformly in $\mathbf{x}$ with
$\|\mathbf{x}\|_2=1$,
    \[
    I_k(\alpha)\geq \inf_{\mathbf{x}\in \mathbb{R}^k:\|\mathbf{x}\|_2=1}
    \sup_{|t|\leq \epsilon}
   \left( t\alpha
    -\log \expec[e^{t S_{\mathbf{x},1}^2}]\right)
    \geq \inf_{\mathbf{x}\in \mathbb{R}^k:\|\mathbf{x}\|_2=1}
    \sup_{|t|\leq \epsilon} \left( t(\alpha-1)
    -\frac{t^2}{2}\min\{2, {\rm Var}(C_{11}^2)\} \right)
    >0,
    \]
the latter bound holding for every $\alpha\neq 1$ when
${\rm Var}(C_{11}^2)>0$. This completes the proof
that $I_k(\alpha)>0$ for all $\alpha\neq 1$ when 
there exists an $\epsilon>0$ such that $\expec[e^{\epsilon C_{11}^2}]<\infty$
and ${\rm Var}(C_{11}^2)>0$.

We continue by proving (\ref{ratemaxa})--(\ref{rateminb}).
The proof for $\lambda_{\rm max}$ is similar to the one for
$\lambda_{\rm min}$, so we will focus on the latter.
To obtain the upper bound of the rate of (\ref{Pmin}),
we use that for any $\mathbf{x'}$ with $\|\mathbf{x'}\|_2=1$,
    \beq
    \mathbb{P}(\lambda_{\rm min} \leq \alpha)=
    \prob(\exists \mathbf{x}: \langle \mathbf{x}, W\mathbf{x}\rangle \leq \alpha)
    \geq \prob(\langle \mathbf{x'}, W\mathbf{x'}\rangle \leq \alpha).
    \label{kieseenx}
    \eeq
Now insert (\ref{Sxrew}). Since $\mathbf{x'}$ is fixed, the $S_{\mathbf{x'},i}^2$ are
i.i.d.\ variables, and we can apply Cram\'er's Theorem to obtain the upper
bound for the rate function for fixed $\mathbf{x'}$. This yields that,
for every $\mathbf{x'}$, we have
    \beq
    -\liminf_{n\rightarrow \infty}
    \frac{1}{n}\log \mathbb{P}(\lambda_{\rm min} \leq \alpha)
    \leq \sup_{t} \left( t\alpha
    -\log \expec[e^{t S_{\mathbf{x'},1}^2}]\right).
    \eeq
If we maximize the right hand side over $\mathbf{x'}$, then
we arrive at $I_k(\alpha)$ as the upper bound, and we have proved
(\ref{ratemina}). The proof for (\ref{ratemaxa}) is identical.

We are left to prove the lower bounds in (\ref{ratemaxb}) and
(\ref{rateminb}). For this, we wish to sum over all possible $\mathbf{x}$.
We approximate the sphere $\|\mathbf{x}\|_2 = 1$ by a finite set of
vectors $\mathbf{x}^\tinbr{j}$ with $\|\mathbf{x}^{\tinbr{j}}\|_2 = 1$,
such that the distance between two of these vectors is at most $d$,
and observe that
    \begin{eqnarray}
    \nonumber |\langle \mathbf{x}, W\mathbf{x}\rangle
    -\langle \mathbf{x}^{\tinbr{j}}, W\mathbf{x}^{\tinbr{j}} \rangle|
    & = &
    |\langle (\mathbf{x} - \mathbf{x}^{\tinbr{j}}), W\mathbf{x}\rangle
    +\langle \mathbf{x}^{\tinbr{j}}, W(\mathbf{x}-\mathbf{x}^{\tinbr{j}}) \rangle|\\
    \nonumber & = &
    |\langle \mathbf{x}, W(\mathbf{x} - \mathbf{x}^{\tinbr{j}})\rangle
    +\langle \mathbf{x}^{\tinbr{j}}, W(\mathbf{x}-\mathbf{x}^{\tinbr{j}}) \rangle| \\
    \nonumber & \leq &
    (\|\mathbf{x}\|+\|\mathbf{x}^{\tinbr{j}}\| ) \cdot\|W\|\cdot\|\mathbf{x}-\mathbf{x}^{\tinbr{j}}\|
    \nonumber \leq  2\lambda_{\rm max} d.
    \end{eqnarray}

We need that $\lambda_{\rm max} \leq \kappa k$, with $\kappa$ some large enough constant,
with sufficiently high probability, which we will prove first.
We have that $\lambda_{\rm max} \leq T_{\sss W}$, where $T_{\sss W}$ is the trace of $W$,
since $W$ is non-negative. Note that
    \beq
    T_{\sss W} =  \frac 1n \sum_{i=1}^n \sum_{m=1}^k C_{mi}^2.
    \eeq
Thus, $T_{\sss W}$ is a sum of $nk$ i.i.d.\ variables.

Since $\expec [e^{t C_{11}^2}]<\infty$
for all $t\leq \epsilon$, we can use Cram\'er's Theorem for $T_{\sss W}$.
Therefore, for any $\kappa$, by the
Chernoff bound,
    \beq
    \prob(T_{\sss W} >\kappa k)\leq e^{-nk I_{C^2}(\kappa)},
    \eeq
where
    \beq
    I_{C^2}(a) =\sup_{t} \left( ta -\log\mathbb{E}[e^{tC_{11}^2}]\right).
    \eeq
Since $\expec[C_{11}^2]={\rm Var}(C_{11})=1$,
we have that $I_{C^2}(\kappa)>0$ for any $\kappa>1$. Therefore, by picking $\kappa>1$
large enough, we can make $k I_{C^2}(\kappa)$ arbitrarily large.
If we take $k I_{C^2}(\kappa)$ larger than $I_k(\alpha-\eps)$,
according to (\ref{largestexponentwins}), this will not influence
the result. (Note that when $I_k(\alpha-\eps)=\infty$ for all $\eps>0$,
then we can also let $k I_{C^2}(\kappa)$ tend to infinity by taking
$\kappa\rightarrow \infty$.)

It follows that
\begin{eqnarray}
  \nonumber \mathbb{P}(\lambda_{\rm min} \leq \alpha) & \leq & \prob(\exists \mathbf{x}^{\tinbr{j}}:
    \langle \mathbf{x}^{\tinbr{j}}, W\mathbf{x}^{\tinbr{j}}\rangle \leq \alpha+2d \kappa k)
    +\prob(T_{\sss W} >\kappa k)\\
  \nonumber & \leq & \sum_j \prob(\langle \mathbf{x}^{\tinbr{j}}, W\mathbf{x}^{\tinbr{j}}\rangle \leq
      \alpha+2d \kappa k)+\prob(T_{\sss W} >\kappa k)\\
 & \leq & N_d \sup_{\mathbf{x}^{\tinbr{j}}} \prob(\langle \mathbf{x}^{\tinbr{j}}, W\mathbf{x}^{\tinbr{j}}\rangle \leq
     \alpha + 2d \kappa k)+\prob(T_{\sss W} >\kappa k),
     \label{metexptermpjes}
\end{eqnarray}
with $N_d$ the number of vectors in the finite approximation of the
sphere. The above bound is valid for every choice of $\kappa, k, \alpha$ and $d$.

We write $\eps=2d \kappa k$ and will later let $\eps\downarrow 0$.
Then, applying the largest-exponent-wins principle
for $\kappa>0$ large enough, as well as Cram\'er's Theorem
together with (\ref{Sxrew}),
we arrive at
    \begin{eqnarray}
    -\limsup_{n\rightarrow \infty}
    \frac 1n \log \mathbb{P}(\lambda_{\rm min} \leq \alpha)
    &\geq& \inf_{\mathbf{x}^{\tinbr{j}}} \sup_t \left( t(\alpha+\eps)
    -\log \expec[e^{t S_{\mathbf{x^{\tinbr{j}}},1}^2}]\right) +
    \liminf_{n\rightarrow \infty} \frac1n \log N_d
    \nonumber\\
    &\geq& I_k(\alpha+\eps) +\liminf_{n\rightarrow \infty} \frac1n \log N_d.
    \end{eqnarray}
In a similar way, we obtain that
    \begin{eqnarray}
    -\limsup_{n\rightarrow \infty}
    \frac 1n \log \mathbb{P}(\lambda_{\rm max} \geq \alpha)
    &\geq& \inf_{\mathbf{x}^{\tinbr{j}}} \sup_t \left( t(\alpha-\eps)
    -\log \expec[e^{t S_{\mathbf{x^{\tinbr{j}}},1}^2}]\right) +
    \liminf_{n\rightarrow \infty} \frac1n \log N_d
    \nonumber\\
    &\geq& I_k(\alpha-\eps) +\liminf_{n\rightarrow \infty} \frac1n \log N_d,
    \end{eqnarray}
where we take $d$ so small that $\alpha - \eps>0$.

A simple overestimation of $N_d$ is obtained
by first taking $[-1,1]^k \subset \mathbb{R}^{k}$
around the origin, and laying a grid on this cube
with grid length $\frac 1L$. We then normalize the centers of these
cubes to have norm 1. The finite set of vectors consists
of the centers of the small cubes of width $2/L$. In this case,
    \beq
    d \leq \frac{3\sqrt{k}}{L}, \hspace{1cm}
    \mbox{and}\hspace{1cm}  N_d \leq L^{k}.
    \label{overestimation}
    \eeq
Indeed, the first bound follows since, for any vector $\mathbf{x}$,
there exists a center of a small cube for which all coordinates are at most
$1/L$ away. Therefore, the distance of to this center is
at most $\frac{\sqrt{k}}{L}$. Since $\mathbf{x}$ has norm 1,
the norm of the center of the cube is in between $1-\frac{\sqrt{k}}{L}$ and
$1+\frac{\sqrt{k}}{L}$, and we obtain that the distance of $\mathbf{x}$
to the normalized center of the small cube is at most
    \beq
    d\leq \frac{\sqrt{k}}{L} +\frac{\frac{\sqrt{k}}{L}}{1-\frac{\sqrt{k}}{L}}
    \leq 3\frac{\sqrt{k}}{L},
    \eeq
when $\frac{\sqrt{k}}{L}\leq 1/2$. For this choice, we have $\eps=6\kappa k^{3/2}/L,$
which we can make small by taking $L$ large.

We conclude that, for any $L<\infty$, $\lim_{n\rightarrow \infty} \frac1n \log N_d=0$, so that,
for any $\kappa>1$ sufficiently large,
    \begin{equation}
    -\limsup_{n\rightarrow \infty} \frac 1n \log \mathbb{P}(\lambda_{\rm min} \leq \alpha)
    \geq I_k(\alpha+\eps),
    \end{equation}
and
    \begin{equation}
    -\limsup_{n\rightarrow \infty} \frac 1n \log \mathbb{P}(\lambda_{\rm max} \geq \alpha)
    \geq I_k(\alpha-\eps),
    \end{equation}
when the respective right-hand sides are finite. Since the above statement
is true for any $\eps$, we can take $\eps\downarrow 0$
by letting $L\uparrow \infty$. When the right-hand side are infinite, then
we conclude that also the left-hand sides can be made arbitrarily
large by letting $L\uparrow \infty$. This completes the proof
of (\ref{ratemaxb}) and
(\ref{rateminb}).
\qed

\vskip0.5cm

\noindent
To see (\ref{extensiontwoa})--(\ref{extensiontwob}),
we follow the above proof. We first note that
the eigenvectors corresponding to $\lambda_{\rm max}$ and
$\lambda_{\rm min}$ are orthogonal. Therefore, we obtain that
    \beq\mathbb{P}(\lambda_{\rm max} \geq \alpha,
    \lambda_{\rm min} \leq \beta)= \prob\big(\exists \mathbf{x,y}: \|\mathbf{x}\|_2 = \|\mathbf{y}\|_2=1,
    \langle \mathbf{x}, \mathbf{y}\rangle=0,
    \langle \mathbf{x}, W\mathbf{x}\rangle \geq \alpha,
    \langle \mathbf{y}, W\mathbf{y}\rangle \leq \beta\big).
    \eeq
We now proceed as above, and for the lower bound pick any $\mathbf{x}, \mathbf{y}$
satisfying the requirements in the probability on the right hand side.
The upper bound is slightly harder. For this, we need to pick a
finite approximation for the choices of $\mathbf{x}$ and $\mathbf{y}$
such that $\|\mathbf{x}\|_2 = \|\mathbf{y}\|_2=1$ and $\langle \mathbf{x},
\mathbf{y}\rangle=0.$ We will now show that we can do this in such a way
that the total number of pairs $\{\mathbf{x}^{\tinbr{i}}, \mathbf{y}^{\tinbr{i,j}}\}_{i,j\geq 1}$
is bounded by $N_d^2$, where $N_d$ is as in (\ref{overestimation}).

We pick $\{\mathbf{x}^{\tinbr{i}}\}_{i\geq 1}$ as in the above proof.
Then, for fixed $\mathbf{x}^{\tinbr{i}}$, we define a finite
number of $\mathbf{y}$ such that $\langle \mathbf{x}^{\tinbr{i}},
\mathbf{y}\rangle=0.$ For this, we consider, for fixed $\mathbf{x}^{\tinbr{i}}$,
only those cubes of width $\frac 1L$ around an
$\mathbf{x}^{\tinbr{j}}$, for some $j$, that contain at least one element
$\mathbf{z}$ having norm 1 and such that $\langle \mathbf{z}, \mathbf{x}^{\tinbr{i}}
\rangle =0$. Fix one of such cubes.
If there are more such $\mathbf{z}$ in this cube around $\mathbf{x}^{\tinbr{j}}$,
then we pick the unique element that is closest to $\mathbf{x}^{\tinbr{j}}$.
We denote this element by $\mathbf{y}^{\tinbr{j,i}}$.
The set of these elements $\mathbf{y}^{\tinbr{j,i}}$ will be
denoted by $\{\mathbf{y}^{\tinbr{j,i}}\}_{i\geq 1}$. The finite subset
of the set $\|\mathbf{x}\|_2 = \|\mathbf{y}\|_2=1$ and $\langle \mathbf{x},
\mathbf{y}\rangle=0$ then consists of $\{\mathbf{x}^{\tinbr{i}},
\mathbf{y}^{\tinbr{i,j}}\}_{i,j\geq 1}$.

We clearly have that every $\mathbf{x}$ and $\mathbf{y}$ with
$\|\mathbf{x}\|_2 = \|\mathbf{y}\|_2=1$ and $\langle \mathbf{x},
\mathbf{y}\rangle=0$ can be approximated by a pair $\mathbf{x}^{\tinbr{j}}$
and $\mathbf{y}^{\tinbr{j,i}}$ such that $\|\mathbf{x}-\mathbf{x}^{\tinbr{j}}\|_2\leq
d$ and $\|\mathbf{y}-\mathbf{y}^{\tinbr{i,j}}\|_2\leq
2d$. Then we can complete the proof as above.

\subsection{Special case: Wishart matrices}
\label{wishart}
To give an example, we go to Wishart matrices, for which $C_{ij}$ are
i.i.d.\ standard normal. In this case, we can compute $I_k(\alpha)$ and
$I_k(\alpha,\beta)$ explicitly. To compute $I_k(\alpha)$, we note that,
for any $\mathbf{x}$ such that $\|\mathbf{x}\|_2=1$, we have that
$S_{\mathbf{x},1}$ is standard normal. Therefore,
    \beq
    \expec[e^{t S_{\mathbf{x},1}^2}]=\frac{1}{\sqrt{1-2t}},
    \eeq
so that
    \beq
    I_k(\alpha) = \sup_{t} \Big(t\alpha-\log\Big(\frac{1}{\sqrt{1-2t}}\Big)\Big).
    \eeq
In order to compute $I_k(\alpha)$, we note that
the maximization problem over $t$ in $\sup_{t} t\alpha
-\log\big(\frac{1}{\sqrt{1-2t}}\big)$
is straightforward, and yields $t^*=\frac 12-\frac {1}{2\alpha}$
and $I_k(\alpha)=\frac 12 (\alpha-1 - \log\alpha)$.
Note that $I_k(\alpha)$ is independent of $k$.
In particular, we see that $\alpha \mapsto I_k(\alpha)$ is
continuous, which leads us to the following corollary:
    \begin{corollary}
    Let $C_{ij}$ be independent standard normals. Then,\\
    (a) for all $\alpha\geq 1$ and fixed $k\geq 2$
    \beq
    \label{ratemaxWM}
    \lim_{n\rightarrow \infty}-\frac 1n \log \mathbb{P}(\lambda_{\rm max}
    \geq \alpha) = \frac 12 (\alpha-1 - \log\alpha),
    \eeq
    (b) for all $0 \leq \alpha \leq 1$ and fixed $k\geq 2$
    \beq
    \label{rateminWM}
    \limsup_{n\rightarrow \infty}
    -\frac 1n \log \mathbb{P}(\lambda_{\rm min} \leq \alpha) =\frac 12 (\alpha-1 - \log\alpha).
    \eeq
    \end{corollary}

We next turn to the computation of $I_k(\alpha,\beta)$.
When $\mathbf{x}$ and $\mathbf{y}$ are
such that $\|\mathbf{x}\|_2=\|\mathbf{y}\|_2=1$ and $\langle \mathbf{x}, \mathbf{y}\rangle=0,$
then we have that $(S_{\mathbf{x},1}, S_{\mathbf{y},1})$ are normally distributed.
It can easily be seen that $\expec[S_{\mathbf{x},1}]=0, \expec[S_{\mathbf{x},1}^2]
=\|\mathbf{x}\|_2^2=1$, so that $S_{\mathbf{x},1}$ and $S_{\mathbf{y},1}$
are standard normal. Moreover, $\expec[S_{\mathbf{x},1}S_{\mathbf{y},1}]
=\langle \mathbf{x}, \mathbf{y}\rangle=0$, so that $(S_{\mathbf{x},1}, S_{\mathbf{y},1})$
are in fact independent standard normal random variables. Therefore,
    \beq
    \expec[e^{t S_{\mathbf{x},1}^2+s S_{\mathbf{y},1}^2}]
    =\expec[e^{t S_{\mathbf{x},1}^2}]\expec[e^{s S_{\mathbf{y},1}^2}]
    =\frac{1}{\sqrt{1-2t}}\frac{1}{\sqrt{1-2s}},
    \eeq
and, for $\alpha \in [0,1]$ and $\beta\geq 1$,
    \beq
    I_k(\alpha,\beta)=\sup_{s,t}
    \left( t\alpha+s\beta-\log\Big(\frac{1}{\sqrt{1-2t}}\Big)-\log\Big(\frac{1}{\sqrt{1-2s}}\Big)\right)
    =I_k(\alpha)+I_k(\beta),
    \eeq
so that the exponential rate of the probability that $\lambda_{\rm max}
\geq \alpha$ and $\lambda_{\rm min} \leq \beta$ is the exponential rate of
the product of the probabilities that $\lambda_{\rm max}
\geq \alpha$ and $\lambda_{\rm min} \leq \beta$.
This remarkable form of independence seems to be true
only for Wishart matrices.

The above considerations lead to the following corollary:
    \begin{corollary}
    Let $C_{ij}$ be independent standard normals. Then,\\
    \beq
    \lim_{n\rightarrow \infty}-\frac 1n \log \mathbb{P}(\lambda_{\rm max}
    \geq \alpha, \lambda_{\rm min} \leq \beta) =
    \frac 12 (\alpha-1 - \log\alpha)
    +\frac 12 (\beta-1 - \log\beta).
    \eeq
    \end{corollary}

In the sequel, we will, among other things, investigate cases where, for $k\rightarrow \infty$,
the rate function $I_k(\alpha)$ for general $C_{ij}$ converges to the Gaussian limit
$I_{\sss \infty}(\alpha)=\frac 12 (\alpha-1 - \log\alpha)$.

\section{Asymptotics for the eigenvalues for symmetric and bounded entries of $C$}
\label{special}
In this section, we investigate the case where $C_{mi}$ is
symmetric around 0 and $|C_{mi}| < M<\infty$ almost surely, or $C_{mi}$ is standard normal.
To emphasize the role of $k$, we will denote the law of $W$ for a given $k$
by $\prob_k$. We define the extension to $k=\infty$ of $I_{k}(\alpha)$ to be
    \beq
    \label{Iinftydef}
    I_{\sss \infty}(\alpha) =
    \inf_{\mathbf{x}\in \ell^2(\mathbb{N}):\|\mathbf{x}\|_2=1} \sup_t \left( t\alpha
    -\log \expec[e^{t S_{\mathbf{x},1}^2}]\right),
    \eeq
where $\ell^2(\mathbb{N})$ is the space of all infinite square-summable
sequences, with norm $\|\mathbf{x}\|_2=\sqrt{\sum_{i=1}^{\infty} \mathbf{x}_i^2}$.
The main result in this section is the following theorem:

\begin{theorem}\label{c}
Suppose that $C_{mi}$ is symmetric around zero and that
$|C_{mi}| < M<\infty$ almost surely, or $C_{mi}$ is standard normal. Then,
for all $k_n \to \infty$ such that $k_n=
o(\frac{n}{\log\log n})$,\\
(a) for all $\alpha\geq 1$,
    \beq
    \label{ratemaxinfa}
    \liminf_{n\rightarrow\infty}-\frac 1n \log \mathbb{P}_{k_n}(\lambda_{\rm max}
    \geq \alpha) \leq I_{\sss \infty}(\alpha),
    \eeq
and
    \beq
    \label{ratemaxinfb}
    \limsup_{n\rightarrow\infty}-\frac 1n \log \mathbb{P}_{k_n}(\lambda_{\rm max}
    \geq \alpha) \geq \lim_{\eps\downarrow 0}I_{\sss \infty}(\alpha -\eps),
    \eeq
(b) for all $0<\alpha\leq 1$,
    \beq
    \label{ratemininfa}
    \liminf_{n\rightarrow\infty}-\frac 1n \log \mathbb{P}_{k_n}(\lambda_{\rm min} \leq \alpha)
    \leq I_{\sss \infty}(\alpha),
    \eeq
and
    \beq
    \label{ratemininfb}
    \liminf_{n\rightarrow\infty}-\frac 1n \log \mathbb{P}_{k_n}(\lambda_{\rm min} \leq \alpha)
    \geq \lim_{\eps\downarrow 0} I_{\sss \infty}(\alpha+\eps).
    \eeq
\end{theorem}

A version of this result has also been published in a conference proceeding
\cite{fey}, for the special case
$C_{mi} = \pm 1$, each with probability $1/2$, and where the restriction
on $k_n$ was $k_n={\cal O}(\frac{n}{\log{n}})$.
Unfortunately, there is a technical error in the proof,
and below we present the corrected proof. In order to do so, we will
rely on explicit lower bounds for $I_k(\alpha)$ for $\alpha\geq 1$.

A priory, it is not obvious that the limit $I_{\sss \infty}(\alpha)$ is strictly
positive for $\alpha\neq 1$. However, in the examples we will investigate
later on, such as $C_{mi} = \pm 1$ with equal probability, we will see
that indeed $I_{\sss \infty}(\alpha)>0$ for $\alpha\neq 1$. Possibly, such a result
can be shown more generally.

The following proposition is
instrumental in the proof of Theorem \ref{c}:

\begin{proposition}\label{Iuitrekenen} Assume that $C_{mi}$ is
symmetric around zero and that $|C_{mi}| < M<\infty$ almost surely, or $C_{mi}$ is standard normal.
Then, for all $k$, $\alpha \geq M^2$ and $\mathbf{x}$ with
$\|\mathbf{x}\|_2=1$,
    \beq
    \mathbb{P}_k(\langle \mathbf{x}, W\mathbf{x}\rangle \geq \alpha)
    \leq e^{-nJ_k(\alpha)},
    \eeq
where
    \beq
    J_k(\alpha)=
    \frac 12 \left(\frac{\alpha}{M^2}-1 - \log\frac{\alpha}{M^2}\right).
    \eeq
\end{proposition}

In the case where $C_{mi} = \pm 1$, for which $M>1$, we will present an
improved version of this bound, valid when
$\alpha \geq 1/2$, in Theorem \ref{Iuitrekenen3} below.

\subsection{Proof of Proposition \ref{Iuitrekenen}}
Throughout this proof, we fix $\mathbf{x}$ with $\|x\|_2=1$.
We use (\ref{Sxrew}) to bound, for every $t\geq 0$
and $k\in \mathbb{N}$, by the Markov inequality,
    \beq
    \mathbb{P}_k(\langle \mathbf{x}, W\mathbf{x}\rangle \geq \alpha)
    = \mathbb{P}_k(e^{t \sum_{i=1}^n S_{\mathbf{x},i}^2}\geq  e^{nt\alpha})
    \leq e^{-n\big(\alpha t-\log\expec_{k_n}[e^{t S_{\mathbf{x},1}^2}]\big)}.
    \eeq
We claim that for all $0\leq t\leq \frac{1}{M^2}$,
    \beq
    \expec_{k_n}[e^{t S_{\mathbf{x},i}^2}]\leq \frac{1}{\sqrt{1-2M^2t}}.
    \label{boundSx}
    \eeq
In the case of Wishart matrices, for which $S_{\mathbf{x},i}$ has a standard
normal distribution, (\ref{boundSx}) holds with equality for $M=1$.

We first note that (\ref{boundSx}) is proven in \cite[Section IV]{hofstad},
for the case that $C_{ij} = \pm 1$ with equal probability. For any $k$ and $\mathbf{x}$,
the bound is even valid for all $-1/2 \leq t\leq 1/2$.
We now extend the case where $C_{ij}=\pm1$ to the case where $C_{ij}$ is symmetric around
zero and satisfies $|C_{ij}|< M$ almost surely.

We write $C_{ij}=A_{ij}C_{ij}^*$, where $A_{ij}=|C_{ij}|< M$ a.s.\ and
$C_{ij}^* = {\rm sign}(C_{ij})$. Moreover, $A_{ij}$ and $C_{ij}^*$ are independent,
since $C_{ij}$ has a symmetric distribution around zero.
Thus, we obtain that $S_{\mathbf{x},i} = S_{A_i \mathbf{x},i}^*$, where
$(A_i \mathbf{x})_j=A_{ij}\mathbf{x}_j$, and
    \beq
    S_{\mathbf{y},i}^*=\sum_{j=1}^k C^*_{ij} \mathbf{y}_j.
    \eeq
For $S_{\mathbf{y},i}^*$ we know that (\ref{boundSx}) is proven.
Therefore,
    \beq
    \expec_k[e^{t S_{A_i \mathbf{x},i}^2}] \leq
    \expec_k\Big[\frac{1}{\sqrt{1-2t \|A_i \mathbf{x}\|_2^2}}
    \Big]
    \eeq
for all $t$ such that $ -1/2 \leq t \|A_i \mathbf{x}\|_2 \leq 1/2$ almost surely.
When $\|\mathbf{x}\|^2_2 = 1$, we have that
    \beq
    0\leq \|A_i \mathbf{x}\|_2 < M \qquad \mbox{a.s.}
    \eeq
Therefore, $\expec_k[e^{t S_{A_i \mathbf{x},i}^2}]
\leq \frac{1}{\sqrt{1-2M^2t\|\mathbf{x}\|_2^2}}$ for all
$0\leq t M^2 \|\mathbf{x}\|_2^2 \leq 1/2$. Thus, we arrive at
    \beq
    \label{Pkbd}
    \mathbb{P}_k(\langle \mathbf{x}, W\mathbf{x}\rangle \geq \alpha)
    \leq e^{-n\big(\sup_{0\leq t\leq 1/M^{2}}
    \left( t\alpha - \log \frac{1}{\sqrt{1-2M^2t}}\right)\big)}.
   \eeq
Note that since $\|\mathbf{x}\|_2=1$, the bound is independent of $\mathbf{x}$.
performing the maximum over $t$ on the right-hand side of (\ref{Pkbd})
over $t$ yields $t^*=\frac{1}{2M^2}-\frac {1}{2\alpha}$, and inserting this
value $t^*$ in the right-hand side gives the result.
\qed

\subsection{Proof of Theorem \ref{c}}
The proof is similar to that of Theorem \ref{prop-LDeig}.
For the proofs of (\ref{ratemaxinfa}) and (\ref{ratemininfa}),
we again use (\ref{kieseenx}), but now choose an
$\mathbf{x'}$ of which only the first $k$ components
are non-zero. This leads to, using that $k_n\rightarrow \infty$,
so that $k_n\geq k$ for $n$ sufficiently large,
    \beq
    \liminf_{n\rightarrow \infty}
    -\frac 1n \log \mathbb{P}_{k_n}(\lambda_{\rm max}
    \geq \alpha) \leq \sup_{t} \left(t\alpha
    -\log \expec[e^{t S_{\mathbf{x'},1}^2}]\right).
    \eeq
Maximizing over all $\mathbf{x'}$ of which only the first $k$ components
are non-zero leads to
    \beq
    \liminf_{n\rightarrow \infty}
    -\frac 1n \log \mathbb{P}_{k_n}(\lambda_{\rm max}
    \geq \alpha) \leq I_k(\alpha),
    \eeq
where this bound is valid for all $k\in \mathbb{N}$. We next claim that
\beq
\lim_{k\rightarrow \infty} I_k(\alpha)=I_{\sss \infty}(\alpha).
\label{knaaroneindig}
\eeq
For this, we first note that the sequence $k\mapsto I_{k}(\alpha)$
is non-increasing and non-negative, so that it has a pointwise limit.
Secondly, $I_{k}(\alpha) \geq I_{\sss \infty}(\alpha)$ for all $k$, since
the possible choices of $\mathbf{x}$ in (\ref{Iinftydef}) is larger than the
one in (\ref{Ikdef}).
Now it is not hard to see that $\lim_{k\rightarrow \infty} I_{k}(\alpha)
=I_{\sss \infty}(\alpha),$ by splitting into the two cases depending on
whether the infimum over $\mathbf{x}$ in (\ref{Iinftydef}) is attained
or not. This completes the proof of (\ref{ratemaxinfa}) and (\ref{ratemininfa}).

For the proof of (\ref{ratemaxinfb}) and (\ref{ratemininfb}),
we adapt the proof (\ref{ratemaxb}) and (\ref{rateminb}).
As in the proof of Theorem \ref{prop-LDeig}(a-b), we wish to
show that the terms $\frac 1n \log{N_d}$ and $2d\lambda_{\rm max}$
vanish when we take the logarithm of (\ref{metexptermpjes}), divide by $n$ and let
$n\rightarrow \infty$. However, this time we wish to let $k_n \rightarrow \infty$
as well, for $k_n$ as large as possible. We will have $k_n=o(n)$ in mind.

The overestimation (\ref{overestimation}) can be improved using an
upper bound for the number $M_{\sss R}=N_{\sss 1/R}$ of spheres of radius $1/R$ needed
to cover the surface of a $k$-dimensional sphere of radius 1
(Rogers (1963)), when $k \to \infty$,
    \beq
    M_{\sss R} = 4 k \sqrt{k} R^k \left(\log{k} +\log\log{k} +
    \log(R)\right)(1+\mathcal{O}(1/\log{k})) \equiv f(k,R)R^k.
    \label{rogersbound}
    \eeq
This bound is valid for $R > \sqrt{\frac{k}{k-1}}$. Since we use
small spheres this time, $d \leq 1/R$.

We can also improve the
upper bound for $\lambda_{\rm max}$. For any $\Omega_n>1$,
which we will choose appropriately later on,
we split
    \begin{eqnarray}
    P_{\rm min}(\alpha)&\leq& \prob(\lambda _{\rm min} \leq \alpha,
    \lambda_{\rm max}\leq \Omega_n) + P_{\rm max}(\Omega_n),\\
    P_{\rm max}(\alpha)&=&\prob(\alpha \leq \lambda_{\rm max}\leq \Omega_n) + P_{\rm max}(\Omega_n).
    \label{lambdagrenzen}
    \end{eqnarray}
We first give a sketch of the proof, omitting the details.
The idea is that the first term of these expressions will yield
the rate function $I_{\sss \infty}(\alpha)$. The term $P_{\rm max}(\Omega_n)$ has
an exponential rate which is $\mathcal{O}(\Omega_n) - \mathcal{O}(\frac{k_n \log k_n}{n})$,
and, since $\frac{k_n \log k_n}{n}=o(\log{n})$, can thus be made arbitrarily
large by taking $\Omega_n=K \log{n}$ with $K>1$ large enough.
This means that we can choose $\Omega_n$ large enough to make this rate
disappear according to the largest-exponent-wins principle
(\ref{largestexponentwins}). We will need different choices
of $R$ for the two terms. We will now give the details of the proof.

We first bound $P_{\rm max}(\Omega_n)$ of (\ref{lambdagrenzen}),
using (\ref{metexptermpjes}). In (\ref{metexptermpjes}), we choose
$\kappa=M^2$. This leads to
    $$
    P_{\rm max}(\Omega_n) \leq M_{\sss R} \sup_\mathbf{x}
    \prob(\langle \mathbf{x},W\mathbf{x} \rangle \geq \Omega_n - 2d M^2 k_n) + \prob(T_{\sss W} \geq M^2 k_n),
    $$
where the supremum over $\mathbf{x}$ runs over the centers
of the small balls. Inserting (\ref{rogersbound}), choosing $R = k_n$ and
using $d\leq 1/R$, this becomes
    $$
    P_{\rm max}(\Omega_n) \leq f(k_n,k_n) k_n^{k_n} \sup_\mathbf{x}
    \prob(\langle \mathbf{x},W\mathbf{x} \rangle \geq \Omega_n - 2M^2)+ \prob(T_{\sss W} \geq M^2 k_n).
    $$
Using Proposition \ref{Iuitrekenen}, we find
    \begin{eqnarray}
    \label{essbd1}
    P_{\rm max}(\Omega_n) & \leq f(k_n,k_n) k_n^{k_n}  e^{- \frac 12 n(\frac{\Omega_n}{M^2} - 3 - \log(\frac{\Omega_n}{M^2} - 2))} + e^{-nk_nI_{C_{11}^2}(M^2)}\nonumber \\
 & =f(k_n,k_n)  e^{k_n\log{k_n}- \frac 12 n\big(\frac{\Omega_n}{M^2} - 3 - \log(\frac{\Omega_n}{M^2}-2)\big)} + e^{-nk_nI_{C_{11}^2}(M^2)}.
    \end{eqnarray}
    
We choose $\Omega_n=K\log{n}$ with $K$ so large that
    $$\frac{k_n \log k_n}{n} < \frac 14
    \big(\frac{\Omega_n}{M^2} - 3 - \log(\frac{\Omega_n}{M^2} - 2)\big).
    $$
Therefore, also using that
    $$
    f(k_n,k_n)=e^{o(n\log{n})},
    $$
we obtain
    \beq
    \label{bdOmegalarge}
    P_{\rm max}(\Omega_n) \leq e^{-\frac{K}{4M^2}n\log{n}(1+o(1))} + e^{-nk_nI_{C_{11}^2}(M^2)}.
    \eeq

Next, we investigate the first term of (\ref{lambdagrenzen}).
In this term, we can use $\Omega_n$ as the upper bound for
$\lambda_{\rm max}$. Therefore, again starting with
(\ref{metexptermpjes}), we obtain that, for any $R$,
    \beq
    \label{essbd}
    -\frac 1n \log \prob(\lambda_{\rm min}\leq \alpha, \lambda_{\rm max}\leq \Omega_n)
    \geq -\frac{1}{n} \big[\log M_{\sss R} +
    \sup_{\mathbf{x}}\log \prob\big(\langle \mathbf{x}, W\mathbf{x}\rangle
    \leq \alpha + 2\Omega_n /R\big)\big].
    \eeq
For $\lambda_{\rm max}$, we get a similar expression.
Inserting (\ref{rogersbound}), we need to choose $R$ again.
This time we wish to choose $R=R_n$ to increase in such a way that
$k_n \log{R_n} = o(n)$ and $\Omega_n=K\log{n}=o(R_n)$. For the latter, we need
that $R_n \gg \log{n}$, so that we can only satisfy the first restriction
when $k_n=o(\frac{n}{\log\log{n}})$.
Then this is sufficient to make the term $2\Omega_n /R_n$ disappear
as $k_n \to \infty$, and to make the term $\frac{1}{n} (\log M_{R_n})
={\cal O}(k_n\log{R_n})$ to be $o(n)$, so that it also disappears. Therefore,
for any $R=R_n$ satisfying the above two restrictions,
    \beq
    \label{lb1}
    -\frac 1n \log \prob(\lambda_{\rm min}\leq \alpha, \lambda_{\rm max}\leq \Omega_n)
    \geq I_{k_n}(\alpha+2\Omega_n/R_n)
    +o(1).
    \eeq
Similarly,
    \beq
    \label{lb1b}
    -\frac 1n \log \prob(\alpha \leq \lambda_{\rm max}\leq \Omega_n)
    \geq I_{k_n}(\alpha-2\Omega_n/R_n)
    +o(1).
    \eeq
Moreover, by the fact that $\Omega_n=K\log{n}=o(R_n)$, we have that
$2\Omega_n/R_n\leq \eps$ for all $n$ large enough. By the monotonicity of
$\alpha \mapsto I_{k_n}(\alpha)$, we then have that
    \beq
    \label{Iknconv}
    I_{k_n}(\alpha+2\Omega_n/R_n)\geq I_{k_n}(\alpha+\eps),
    \qquad
    I_{k_n}(\alpha-2\Omega_n/R_n)\geq I_{k_n}(\alpha-\eps).
    \eeq
%
%

Since $\lim_{k\rightarrow \infty} I_{k}(\alpha)=I_{\sss \infty}(\alpha)$ (see (\ref{knaaroneindig})),
putting (\ref{lb1}), (\ref{Iknconv}) and (\ref{bdOmegalarge}) together
and applying the largest-exponent-wins principle (\ref{largestexponentwins}), we
see that the proof follows when
    \beq
    I_{\sss \infty}(\alpha\pm \eps) < \min\{\log{n}\frac{K}{4M^2},k_nI_{C_{11}^2}(M^2)\}.
    \label{voorwaarde}
    \eeq
Both terms are increasing in $n$, as long as $I_{C_{11}^2}(M^2) >0$. This is true for the $C_{11}$ we consider: if $C_{11}$ is
symmetric around zero such that $|C_{11}| < M<\infty$ almost surely, then $I_{C_{11}^2}(M^2) = \infty$, and if $C_{mi}$ is standard normal then $I_{C_{11}^2}(M^2) >0$.

Therefore, (\ref{voorwaarde}) is true when $n$ is large enough, when $I_{\sss \infty}(\alpha\pm \eps)
<\infty$. On the other hand, when $I_{\sss \infty}(\alpha\pm \eps)=\infty$, then we obtain
that the exponential rates converge to infinity, as stated in
(\ref{ratemaxinfb}) and (\ref{ratemininfb}). We conclude that,
for every $\eps>0$,
    \beq
    \label{ratemininfaeps}
    \liminf_{n\rightarrow\infty}-\frac 1n \log \mathbb{P}_{k_n}(\lambda_{\rm min}
    \leq \alpha) \geq I_{\sss \infty}(\alpha+\eps),
    \eeq
and
    \beq
    \label{ratemaxinfaeps}
    \liminf_{n\rightarrow\infty}-\frac 1n \log \mathbb{P}_{k_n}(\lambda_{\rm max}
    \geq \alpha) \geq I_{\sss \infty}(\alpha-\eps),
    \eeq
and letting $\eps\downarrow 0$ completes the proof
for every sequence $k_n$ such that $k_n=o(\frac{n}{\log\log{n}})$.
The proof for $\lambda_{\rm max}$ is identical to the above proof.
\qed

We believe that the above argument can be extended somewhat further, by making
a further split into $K'\log\log{n}\leq \lambda_{\rm max}\leq K\log{n}$ and
$\lambda_{\rm max}\leq K'\log\log{n}$, but we refrain from writing this down.

\subsection{The limiting rate for $k$ large}
In this section, we investigate what happens when we take $k$ large.
In certain cases, we can show that the rate function, which depends on $k$,
converges to the rate function for Wishart matrices. This will be formulated in the
following theorem.

\begin{theorem}\label{Iuitrekenen2}
Assume that $C_{ij}$ satisfies (\ref{Cass}), and, moreover, that
$\phi_{\sss C}(t) \leq e^{t^2/2}$ for all $t$.  Then, for all $\alpha \geq 1$,
and all $k\geq 2$,
    \beq
    \label{bdrate}
    I_k(\alpha) \geq
    \frac 12 (\alpha-1 - \log\alpha),
    \eeq
and, for all $\alpha\geq 1$,
    \beq
    \label{Ikconv}
    \lim_{k\rightarrow \infty}I_k(\alpha)=I_{\sss \infty}(\alpha)=\frac 12 (\alpha-1 - \log\alpha).
    \eeq
Finally, for all $k_n\rightarrow \infty$ such that $k_n=
o(\frac{n}{\log\log n})$ and $\alpha\geq 1$,
    \beq
    \label{LEratebd}
    \lim_{n\rightarrow \infty}-\frac 1n \log \mathbb{P}_{k_n}(\lambda_{\rm max}
    \geq \alpha) = \frac 12 (\alpha-1 - \log\alpha).
    \eeq
\end{theorem}
Note that, in particular, Theorem \ref{Iuitrekenen2} implies that
$I_{\sss \infty}(\alpha)=\frac 12 (\alpha-1 - \log\alpha)>0$ for all
$\alpha>1$.

Theorem \ref{Iuitrekenen2} is a kind of universality result, and shows that, for $k$ large,
the rate functions of certain sample covariance matrices converges to the rate
function for Wishart matrices. An example where $\phi_{\sss C}(t) \leq e^{t^2/2}$ holds
is when $C_{ij}=\pm 1$ with equal probability. We will call this example the Bernoulli case.
A second example is uniform random variable on
$[-\sqrt{3}, \sqrt{3}]$, for which also the variance equals 1. We will prove these bounds below.

Of course, the relation that $\phi_{\sss C}(t) \leq e^{t^2/2}$ for random
variables with mean 0 and variance 1, is equivalent to the statement
that $\phi_{\sss C}(t) \leq e^{t^2\sigma^2/2}$ for a random variable $C$ with
mean $0$ and variance $\sigma^2$. Thus, we will check the condition for
uniform random variables on $[-1,1]$ and for the Bernoulli case.
We will denote the moment generating functions by $\phi_U$ and $\phi_B$.
We start with the second, for which we have that
    \beq
    \phi_B(t) = \cosh(t) =\sum_{n=0}^{\infty} \frac{t^{2n}}{(2n)!}\leq \sum_{n=0}^{\infty} \frac{t^{2n}}{2^n n!}=e^{t^2/2},
    \eeq
since $(2n)! \geq 2^n n!$ for all $n\geq 0$.
The proof for $\phi_U$ is similar. Indeed,
    \beq
    \phi_U(t) = \frac{\sinh(t)}{t} = \sum_{n=0}^{\infty} \frac{t^{2n}}{(2n+1)!}\leq \sum_{n=0}^{\infty} \frac{t^{2n}}{6^n n!}=e^{t^2/6}=e^{t^2\sigma^2/2},
    \eeq
since now $(2n+1)! \geq 6^n n!$ for all $n\geq 0$.
\vskip0.3cm

\noindent
{\bf Proof of Theorem \ref{Iuitrekenen2}.}
Using Theorem \ref{prop-LDeig} and \ref{c}, we claim that it suffices to
prove that, uniformly for $\mathbf{x}$ with $\|\mathbf{x}\|_2=1$
and $t<1$,
    \beq
    \label{bdMGF}
    \expec[e^{t S_{\mathbf{x},i}^2}] \leq \frac{1}{\sqrt{1-2t}}.
    \eeq
We will prove (\ref{bdMGF}) below, and first prove Theorem \ref{Iuitrekenen2}
assuming that (\ref{bdMGF}) holds.

When (\ref{bdMGF}) holds, then, using Theorem \ref{c} and (\ref{Ikdef}),
it immediately follows that (\ref{bdrate}) holds. Here we also use that
$\alpha\mapsto \frac 12 (\alpha-1 - \log\alpha)$ is continuous,
so that the limit over $\eps\downarrow 0$ can be computed.

To prove (\ref{Ikconv}), we take $\mathbf{x}=\frac{1}{\sqrt{k}}(1, \ldots, 1)$,
to obtain that, with $S_k=\sum_{i=1}^k C_{i1}$,
    \beq
    I_k(\alpha) \leq \sup_t \left( t\alpha
    -\log \expec[e^{\frac{t}{k} S_{k}^2}]\right).
    \eeq
We claim that, when $k\rightarrow \infty$, for all $0\leq t<1$,
    \beq
    \label{convMGF}
    \expec[e^{\frac{t}{k} S_{k}^2}] \rightarrow \expec[e^{t Z^2}]=\frac{1}{\sqrt{1-2t}},
    \eeq
where $Z$ is a standard normal random variable. This implies the lower bound
for $I_k(\alpha)$, and thus
(\ref{Ikconv}). Equation (\ref{LEratebd}) follows in a similar way,
also using that $\alpha\mapsto \frac 12 (\alpha-1 - \log\alpha)$ is continuous.

We complete the proof by showing that (\ref{bdMGF}) and (\ref{convMGF}) hold.
We start with (\ref{bdMGF}). We rewrite, for $t\geq 0$,
and writing $Z$ for a standard normal random variable,
    \beq
    \expec[e^{t S_{\mathbf{x},i}^2}] =
    \expec[e^{\sqrt{2t}Z S_{\mathbf{x},i}}]
    =\expec\big[\prod_{j=1}^k \phi_{\sss C}(\sqrt{2t}Z \mathbf{x}_j)\big].
    \eeq
We now use that $\phi_{\sss C}(t) \leq e^{t^2/2}$ to arrive at
    \beq
    \expec[e^{t S_{\mathbf{x},i}^2}] \leq
    \expec\big[\prod_{j=1}^k e^{tZ^2 \mathbf{x}_j^2}\big]
    =\expec[e^{tZ^2}]=\frac{1}{\sqrt{1-2t}}.
    \eeq
This completes the proof of (\ref{bdMGF}). We proceed with (\ref{convMGF}).
We use
    \beq
    \expec[e^{t S_k^2}]=
    \expec\big[\prod_{j=1}^k \phi_{\sss C}(\sqrt{\frac{2t}{k}}Z)\big].
    \eeq
We will use dominated convergence.
By the assumption, we have that $\phi_{\sss C}(\sqrt{\frac{2t}{k}}Z)\leq e^{\frac t{k} Z^2},$
so that $\prod_{j=1}^k \phi_{\sss C}(\sqrt{\frac{2t}{k}}Z)\leq e^{t Z^2}$,
which has a finite expectation when $t<1/2$.
Moreover, $\prod_{j=1}^k \phi_{\sss C}(\sqrt{\frac{2t}{k}}z)$ converges to $e^{t z^2}$
pointwise in $z$. Therefore, dominated convergence proves the claim in (\ref{convMGF}),
and completes the proof.
\qed

\section{The smallest eigenvalue for $C_{ij}=\pm1$}
\label{sec-pm1}

Unfortunately, we are not able to prove a similar result as in Theorem
\ref{Iuitrekenen2} for the smallest eigenvalue. In fact, as we will comment on in
more detail in Section \ref{conjecture} below, we expect the result to be {\it false}
for the smallest eigenvalue, in particular when $\alpha$ is small.
There is one example where we can prove
a partial convergence result, and that is when $C_{ij}=\pm 1$ with equal
probability. Indeed, in this case it is shown in
\cite[Section IV]{hofstad} that (\ref{bdMGF}) holds for all $t\geq -1$. This
leads to the following result, which also implies that 
$I_k(\alpha)>0$ for $\alpha\neq 1$ in the case where ${\rm Var}(C_{11}^2)=0$
(recall also Theorem \ref{prop-LDeig}):

\begin{theorem}\label{Iuitrekenen3}
Assume that $C_{ij}=\pm 1$ with equal probability.
Then, for all $\alpha \geq 1/2$, and all $k\geq 2$,
    \beq
    \label{bdrate2}
    I_k(\alpha) \geq
    \frac 12 (\alpha-1 - \log\alpha),
    \eeq
and, for all $\alpha \geq 1/2$,
    \beq
    \label{Ikconv2}
    I_{\sss \infty}(\alpha)=\frac 12 (\alpha-1 - \log\alpha).
    \eeq
Finally, for all $0 < \alpha \leq 1/2$,
    \beq
    \label{complbd}
    I_k(\alpha) \geq \frac 12 (-\alpha + \log 2).
    \eeq
\end{theorem}

\proof The proof of (\ref{bdrate2}--\ref{Ikconv2}) is identical
to the proof of Theorem \ref{Iuitrekenen2}, now using that (\ref{bdMGF})
holds for all $t\geq -1$. Equation (\ref{complbd}) follows
since $I_k(\alpha)\geq \inf_{\mathbf{x}: \|\mathbf{x}\|_2=1}\left(-\frac{\alpha}{2}
-\log \expec[e^{-\frac{1}{2} S_{\mathbf{x},i}^2}]\right)$ and the bound
on the moment generating function for $t=-1$.
\qed

\subsection{Rate for the probability of one or more zero eigenvalues for $C_{ij}=\pm 1$}
In the above computations, we obtain no control over the probability of
a large deviation of the smallest eigenvalue $\lambda_{\rm min}$.
In this and the next section, we investigate this problem in the case where
$C_{ij}=\pm 1$.
\begin{proposition}\label{prop-lambdanul}
Suppose that $C_{ij}=\pm 1$ with equal probability.
Then, for all $0 < l \leq k-1$, and any $k_n={\cal O}(n^b)$ for some $b$,
\beq
\label{lzeroev}
\lim_{n\rightarrow \infty}-\frac 1n \log \prob_{k_n}(\lambda_1=\ldots=\lambda_l =
0) = l\log{2},
\eeq
where the $\lambda_i$ denote eigenvalues of $W$ arranged in increasing order.
     \end{proposition}

\noindent
\proof
The upper bound in (\ref{lzeroev}) is simple, since, to have $l$
eigenvalues equal to zero, we can take the first $l+1$ columns of $C$
to be equal. For eigenvectors $\mathbf{w}$ of $W$,
    $$ \langle \mathbf{w}, W\mathbf{w}\rangle = \frac 1n \langle
    \mathbf{w}, C C^T \mathbf{w}\rangle = \frac 1n \|C^T
    \mathbf{w}\|_2^2,
    $$
we obtain that $\mathbf{w}$ is an eigenvector with eigenvalue zero
precisely when $\|C^T
\mathbf{w}\|_2 = 0$. When the first $l+1$ columns of $C$ are equal,
then there are $l$ linearly independent vectors for which $\|C^T
\mathbf{w}\|_2 = 0$, so that the multiplicity of the eigenvalue zero
is at least $l$. Moreover, the probability that the first $l+1$ columns
of $C$ are equal is equal to $2^{nl}$.

We prove the lower bound in (\ref{lzeroev}) by induction by showing
that
\beq
\label{lzeroIH}
\lim_{n\rightarrow \infty}-\frac 1n \log \prob_{k_n}(\lambda_1=\ldots=\lambda_l =
0) \geq  l\log{2}.
\eeq
When $l=0$, then the claim is trivial. It suffices to advance the
induction hypothesis.

Suppose that there are $l$ linear independent eigenvectors
with eigenvalue zero. Since the
eigenvectors can be chosen to be orthogonal, it is possible to make linear
combinations, such that the first $l-1$ all have one zero
coordinate $j$, whereas the $l^{\rm th}$ has all
coordinates zero except coordinate $j$. This means that the first $l-1$ eigenvectors
fix some part of $C^T$, but not the $j^{\rm th}$ column. The
$l^{\rm th}$ eigenvector however fixes precisely this column.
Fixing one column of $C^T$ has probability $2^{-n}$. The number of possible
rows $j$ is bounded by $k$, which is turn is bounded by $n^b=e^{o(n)}$. Therefore, we
have
    \beq
    \lim_{n\rightarrow \infty}-\frac 1n \log
    \prob_{k_n}(\lambda_1=\ldots=\lambda_l =
    0) \geq \log{2} + \lim_{n\rightarrow \infty}-\frac 1n \log
    \prob_{k_n}(\lambda_1=\ldots=\lambda_{l-1} =
    \label{inductiestap}
    0).
    \eeq
The claim follows from the induction hypothesis.
\qed

Note that Proposition \ref{prop-lambdanul} shows that (\ref{Ikconv2}) cannot be extended to
$\alpha=0$. Therefore, a changeover takes place between $\alpha=0$ and $\alpha\geq \frac 12$,
where for $\alpha\geq \frac 12$, the rate function equals the one for Wishart matrices, while
for $\alpha=0$, this is not the case. We will comment more on this in Conjecture \ref{conj-small}
below.

\subsection{A conjecture about the smallest eigenvalue for $C_{ij}=\pm 1$}
\label{conjecture}

We have already shown that (\ref{bdrate2}) is sharp.
By Proposition \ref{prop-lambdanul}, (\ref{complbd}) is not sharp, since
$I_{\rm min}(0) =  \log 2$, whereas (\ref{complbd}) only yields
$\lim_{\alpha\downarrow 0} I_{\rm min}(\alpha) \geq  \frac 12 \log 2$.

We can use (\ref{kieseenx}) again with $\mathbf{x'}=\frac{1}{\sqrt{2}}
(1,1,0,\cdots)$. For this vector, $\expec[e^{t S^2_{\mathbf{x'},i}}]
= \frac12(e^{2t}+1$), and calculating the according rate function gives
$I_k(\alpha) \leq I^{\scriptscriptstyle(2)}(\alpha)=
\frac{\alpha}{2} \log \alpha +\frac{2-\alpha}{2} \log(2-\alpha)$,
which implies that
$\lim_{\alpha\downarrow 0} I_k(\alpha) \leq  \log 2$.

It appears that below a certain $\alpha = \alpha^*_k$, the optimal strategy changes
from $\mathbf{x}^{{\scriptscriptstyle (k)}}=\frac{1}{\sqrt{k}} (1,1,\cdots)$
to $\mathbf{x}^{{\scriptscriptstyle(2)}}
=\frac{1}{\sqrt{2}} (1,1,0,\cdots)$. In words, that means that for not too
small eigenvalues of $W$, all entries of $C$ contribute equally to create a
small eigenvalue. However, smaller values of the smallest eigenvalues of $W$
are created by only two columns of $C$, whereas the others are
close to orthogonal. Thus, a change in strategy occurs, which gives rise to a phase
transition in the asymptotic exponential rate for the smallest eigenvalue.

We have the following conjecture:
  \begin{conjecture}
  \label{conj-small}
    For each $k$ and all $\alpha \geq 0$, there exists an $\alpha^*=\alpha^*_k> 0$ so that
    $$
    I_k(\alpha) = I^{\scriptscriptstyle(2)}(\alpha) = \frac{\alpha}{2}\log{\alpha} +
    \frac{(2-\alpha)}{2}\log{(2-\alpha)},
    $$
    for $\alpha \leq \alpha^*_k$. For $\alpha > \alpha^*_k$,
   $$
    I^{\scriptscriptstyle(2)}(\alpha) > I_k(\alpha) \geq I_{\scriptscriptstyle\infty}(\alpha)
    = \frac 12 (\alpha -1 -\log{\alpha}).
    $$
    For $k \to \infty$, the last inequality will become an equality.
    Consequently, $\lim_{k\rightarrow \infty} \alpha^*_k =\alpha^*$, which is
    the positive solution of $I_{\scriptscriptstyle\infty}(\alpha)=I^{\scriptscriptstyle(2)}(\alpha)$.
    \end{conjecture}

For $k = 2$, the conjecture is trivially true, since the two optimal strategies are the same, and the only possible. Note that in the proof of Proposition \ref{prop-lambdanul}, we used that to have a zero eigenvalue, we need two vectors of $C$ to be equal. Thus, the conjecture is also proven for $\alpha = 0$.
Furthermore, with Theorem \ref{Iuitrekenen} and (\ref{boundSx}),
the conjecture follows for all $k$ for $\alpha\geq \frac 12$. We
lack a proof for $0<\alpha<\frac 12$.
Numerical evaluation gives that $\alpha^*_3 \approx 0.425$, and for $k \to \infty$,
$\alpha^*_k \approx 0.253$. We have some evidence that suggests that
$\alpha^*_k$ decreases with $k$.

\newcommand{\smallsup}[1]{{\scriptscriptstyle{(#1)}}}
\section{An application: Mobile Communication Systems}
\label{cdma}
Our results on the eigenvalues of sample covariance matrices was triggered
by a problem in mobile communication systems. In this case, we take the
$C$ matrix as a coding sequence, for which we can assume that the elements
are $\pm 1$. Thus, all our results apply to this case.
In this section, we will describe the consequences of our results
on this problem.

\subsection{Soft-Decision Parallel Interference Cancellation}
In Code Division Multiple Access (CDMA) communication systems,
each of $k$ users multiplies his data signal by an individual
coding sequence. The base station can distinguish the
different messages by taking the inner product of the total signal
with each coding sequence. This is called Matched Filter (MF)
decoding. An important application is mobile telephony. Since, due
to synchronisation problems, it is unfeasible to implement
completely orthogonal codes for mobile users, the decoded messages
will suffer from Multiple Access Interference (MAI). In practice,
pseudo-random codes are used. Designers of decoding schemes are
interested in the probability that a decoding error is made.

In the following, we explain a specific method to iteratively
estimate and subtract the MAI, namely, Soft Decision Parallel Interference
Cancellation (SD-PIC). For more background on SD-PIC, see
\cite{buehrer}, \cite{grant}, \cite{rasmussen1} and \cite{klok},
as well as the references therein. Because this procedure is linear, it can be
expressed in matrix notation. We will show that the possibility of
a decoding error is related to a large deviation of the maximum or
minimum eigenvalue of the code correlation matrix.

To highlight the aspects that are relevant in this article, we
suppose that each user sends only one data bit $b_m\in \{+1, -1\}$, and we omit
noise from additional sources. We can denote all sent data multiplied by
their amplitude in a column vector $\mathbf{Z}$, i.e., $\mathbf{Z}_m=\sqrt{P_m}b_m$,
where $P_m$ is the power of the $m^{\rm th}$ user. The $k$ codes are modeled as
the different rows of length $n$ of the code matrix $C$, consisting of
i.i.d.\ random bits with distribution
    $$
    \mathbb{P}(C_{mi}=+1) =
    \mathbb{P}(C_{mi}=-1) = \frac12.
    $$
Thus, $k$ plays the role of the number of users, while $n$ is the length
of the different codes.

The base station then receives a
total signal  ${\mathbf s} = C^T \mathbf{Z}$. Decoding for user
$m$ is done by taking the inner product with the code of the $m^{\rm th}$ user
$(C_{1m}, \ldots, C_{nm})$, and dividing by $n$. This yields an estimate
$\hat{\mathbf{Z}}_{m}^{\smallsup{1}}$ for the sent signal $\mathbf{Z}_m$. In matrix notation, the
vector $\mathbf{Z}$ is estimated by
    $$\mathbf{\hat{Z}}^{\smallsup{1}} =
    \frac{1}{n} C \mathbf{s} = W \mathbf{Z}.
    $$
Thus, we see that multiplying with the
matrix $W$ is equivalent to the MF decoding scheme.
In order to estimate the signal, we must find
the inverse matrix $W^{-1}$.
From $\mathbf{\hat{Z}}^{\smallsup{1}}$, we estimate the sent bit $b_m$ by
    \beq
    \hat b_m^{\smallsup{1}} = {\rm sign}(\mathbf{\hat{Z}}^{\smallsup{1}}_m)
    \eeq
(where, when $\mathbf{\hat{Z}}^{\smallsup{1}}_m=0$, we toss an independent
fair coin to decide what the value of ${\rm sign}(\mathbf{\hat{Z}}^{\smallsup{1}}_m)$
is). Below, we explain the role of the eigenvalues of $W$ in the
more advanced SD-PIC decoding scheme.

The MF estimate contains MAI. When we write $\mathbf{\hat{Z}} =
\mathbf{Z} + (W-I)\mathbf{Z}$, it is clear that the estimated bit
vector is a sum of the correct bit vector and MAI. In SD-PIC, the
second term is subtracted, with $\mathbf{Z}$ replaced by
$\mathbf{\hat{Z}}$. In the case of multistage PIC, each new
estimate is used in the next PIC iteration. We will now write the
multistage SD-PIC procedure in matrix notation. We number the
successive SD estimates for $\mathbf{Z}$ with an index $s$, where
$s = 1$ corresponds to the MF decoding. In each new iteration, the
latest guess for the MAI is subtracted. The iteration in a
recursive form is therefore:
    \beq
    \mathbf{\hat{Z}}^{\tinbr{s}} = \mathbf{\hat{Z}}^{\tinbr{1}}-
    (W - I)\mathbf{\hat{Z}}^{\tinbr{s-1}}. \label{iteration}
    \eeq
This can be worked out to
    \beq \mathbf{\hat{Z}}^{\tinbr{s}} =
    \sum_{\varsigma=0}^{s-1} (I-W)^{\varsigma}W\mathbf{Z}.
    \label{s-stage}
    \eeq
We then estimate $b_m$ by
    \beq
    \hat b_m^{\smallsup{s}} = {\rm sign}(\mathbf{\hat{Z}}^{\smallsup{s}}_m).
    \eeq

When $s\rightarrow\infty$, the
series $\sum_{\varsigma=0}^{s-1} (I-W)^{\varsigma}$ converges to
$W^{-1}$, as long as the eigenvalues of $W$ are between 0 and 2.
Otherwise, a decoding error is made. This is the crux to
the method, see also \cite{rasmussen2} for the above
matrix computations. When $k = o(n/\log\log n)$, the values
$\lambda_{\rm min} = 0$ and $\lambda_{\rm max} \geq 2$ are large
deviations, and therefore our derived rate functions provide
information on the error probability. In the next section, we will
describe these results, and we will also obtain bounds on the exponential
rate of a bit error in the case that $s$ is fixed and $k$ is large.
For an extensive introduction to CDMA and PIC procedures, we refer
to \cite{klok}.

\subsection{Results for Soft-Decision Parallel Interference Cancelation}
There are two cases that need to be distinguished, namely, the case where
$s\rightarrow \infty$, and the case where $s$ is fixed. We start with the former,
which is simplest. As explained in the previous section, due to the absence of noise,
there can only be bit-errors when $\lambda_{\rm min}=0$ or when $\lambda_{\rm max}\geq 2$.
By (\ref{complbd}), the rate of $\lambda_{\rm min}=0$ is at least $\frac 12 \log{2}\approx 0.35...$,
whereas the rate of $\lambda_{\rm max}\geq 2$ is bounded below by $\frac 12 - \frac 12\log2\approx
0.15...$ The latter bound is weaker, and thus, by the largest-exponent-wins principle,
we obtain the following result:

\begin{theorem}[Bit-error rate for optimal SD-PIC]
\label{thm-OSD}
For all $k$ fixed, or for
$k=k_n\rightarrow \infty$ such that $k_n=
o(\frac{n}{\log\log n}),$
    \beq
    -\frac 1n \log \prob_k\big(\exists m=1, \ldots, k \mbox{ for which }\lim_{s\rightarrow \infty}
    \hat{b}_m^{\tinbr{s}}\neq b_m\big) \geq \frac 12 - \frac 12\log2.
    \eeq
\end{theorem}
We emphasize that in the statement of the result, we write that
$\lim_{s\rightarrow \infty} \hat{b}_m^{\tinbr{s}}\neq
b_m\forall m=1, \ldots, k$ for the statement that either
$\lim_{s\rightarrow \infty} \hat{b}_m^{\tinbr{s}}$
does not exist, or that $\lim_{s\rightarrow \infty} \hat{b}_m^{\tinbr{s}}$
exists, but is unequal to $b_m$. We observe that when $\lambda_{\rm max}>2,$
then
    \beq
    \hat{\mathbf{Z}}^{\tinbr{s}} = \sum_{\varsigma=0}^{s-1} (I-W)^{\varsigma}W\mathbf{Z}
    \eeq
oscillates, so that we can expect there to be errors in every stage.
This is sometimes called the {\it ping-pong effect}
(see \cite{rasmussen2}). Thus, one would expect that
    \[
    -\frac 1n \log \prob_{k_n}\big(\exists m=1, \ldots, k \mbox{ for which }
    \lim_{s\rightarrow \infty} \hat{b}_m^{\tinbr{s}}\neq
    b_m\big)=\frac 12 - \frac 12\log2.
    \]
However, this depends also on the relation between $\mathbf{Z}$ and the eigenvector corresponding
to $\lambda_{\rm max}$. Indeed, when $\mathbf{Z}$ is orthogonal to the eigenvector corresponding
to $\lambda_{\rm max}$, then the equality does not follow. To avoid such problems, we stick to lower bounds
on the rates in this section, rather than asymptotics.

We next go to the case where $s$ is fixed. We again consider the case where
$k$ is large and fixed, or that $k=k_n\rightarrow \infty$. In this case, it can be expected
that the rate converges to 0 as $k\rightarrow \infty$. We already know that
the probability that $\lambda_{\rm max}\geq 2$ or $\lambda_{\rm min}=0$
is exponentially small with {\it fixed} strictly positive lower bound
on the exponential rate. Thus,
we shall assume that $0<\lambda_{\rm min}\leq \lambda_{\rm max}< 2.$
We can then rewrite
    \beq
    \hat{\mathbf{Z}}^{\tinbr{s}} = \sum_{\varsigma=0}^{s-1} (I-W)^{\varsigma}W\mathbf{Z}
    =\big[I- (I-W)^{s}\big]\mathbf{Z}.
    \eeq
For simplicity, we will first assume that $\mathbf{Z}_i=\pm 1$ for all $i=1, \ldots, k$,
which is equivalent to assuming that all powers are equal.
When $s$ is fixed, we cannot have any bit-errors when
    \beq
    \Big|\big((I-W)^{s} \mathbf{Z}\big)_i\Big|<1.
    \eeq
We can bound
    \beq
    \label{restrvep}
    \Big|\big((I-W)^{s} \mathbf{Z}\big)_i\Big| \leq \varepsilon_k^s \|\mathbf{b}\|_2,
    \eeq
where $\varepsilon_k=\max\{1-\lambda_{\rm min},\lambda_{\rm max}-1\}$. Since
$\|\mathbf{b}\|_2=\sqrt{k}$, we obtain that there cannot be any bit-errors when
$\varepsilon_k^s \sqrt{k}<1$. This gives an explicit relation between the bit-errors
and the eigenvalues of a random sample covariance matrix.
By applying the results from the previous two sections,
we obtain the following theorem:

\begin{theorem}[Bit-error rate for finite-stage SD-PIC and $k$ fixed]
\label{thm-fsSDPIC}
For all $k$ such that $k>2^{2s}$,
    \beq
    -\liminf_{n\rightarrow \infty}\frac {1}{n} \log \prob_k\big(\exists m=1, \ldots, k \mbox{ for which }
    \lim_{s\rightarrow \infty} \hat{b}_m^{\tinbr{s}}\neq
    b_m\big) \geq \frac 1{4\sqrt[s]{k}}\Big(1+{\cal O}\big(\frac{1}{\sqrt[s]{k}}\big)\Big).
    \eeq
\end{theorem}

When the signals are different, related results can be obtain in terms
of the minimal and maximal element of $\mathbf{Z}$. We will not write this
case out.
\vskip0.5cm

\proof By the computation in (\ref{restrvep}), there can be no bit-errors when
$1-\lambda_{\rm min}$ and $\lambda_{\rm max}-1$ are both at most $1/\sqrt[2s]{k}$.
Thus,
    \beq
    \label{BEPbdk}
    \prob_k\big(\exists m=1, \ldots, k \mbox{ for which }\lim_{s\rightarrow \infty} \hat{b}_m^{\tinbr{s}}\neq
    b_m\big) \leq \prob_k\big(\lambda_{\rm min}\leq 1-\frac{1}{\sqrt[2s]{k}}\big)
    +\prob_k\big(\lambda_{\rm max}\geq 1+\frac{1}{\sqrt[2s]{k}}\big).
    \eeq
Each of these terms is bounded by, using Theorem \ref{prop-LDeig},
    \beq
    e^{-n \min \Big\{\lim_{\eps\downarrow 0} I_k\big(1-\frac{2}{\sqrt[2s]{k}}-\eps\big),
    I_k\big(1+\frac{2}{\sqrt[2s]{k}}+\eps\big)\Big\}(1+o(1))}.
    \eeq
Since, by Theorem \ref{Iuitrekenen3} and $\alpha\geq \frac 12$, we have that $I_k(\alpha)
\geq I_{\sss \infty}(\alpha)=\frac 12 (\alpha-1 - \log\alpha)$, and
    \beq
    \label{IinftyTaylor}
    I_{\sss \infty}(\alpha) = \frac 14 (\alpha-1)^2 +{\cal O}(|\alpha-1|^3),
    \eeq
the result follows when $k$ is so large that $1-\frac{1}{\sqrt[2s]{k}}>\frac 12$.
The latter is equivalent to $k>2^{2s}$.
\qed

We finally state a result that applied to $k=k_n$:
\begin{theorem}[Bit-error rate for finite-stage SD-PIC and $k=k_n$]
\label{thm-fsSDPICn}
For $k_n = o\big(\frac{n^{\frac{s}{s+1}}}{\log{n}}\big)$,
    \beq
    -\frac {\sqrt[s]{k_n}}{n} \log \prob_{k_n}\big(\exists m=1, \ldots, k \mbox{ for which }
    \lim_{s\rightarrow \infty} \hat{b}_m^{\tinbr{s}}\neq
    b_m\big) \geq \frac 1{4}+{\cal O}(\frac{1}{\sqrt[s]{k_n}}).
    \eeq
\end{theorem}

\proof We use (\ref{BEPbdk}), to conclude that we need
to derive bounds for $\prob_{k_n}\big(\lambda_{\rm min}\leq 1-\frac{1}{\sqrt[2s]{k_n}}\big)$ and
\linebreak
$\prob_{k_n}\big(\lambda_{\rm max}\geq 1+\frac{1}{\sqrt[2s]{k_n}}\big)$. Unfortunately,
the bounds $1-\frac{1}{\sqrt[2s]{k_n}}$ and $1+\frac{1}{\sqrt[2s]{k_n}}$ on the smallest
and largest eigenvalues depend on $n$, rather than being fixed. Therefore, we need
to adapt the proof of Theorem \ref{c}.

We note that, by Theorem \ref{c},
    \[
    \prob_{k_n}(\lambda_{\rm max}\geq 2)
    =e^{-(\frac 12-\frac 12 \log{2})n(1+o(1))}.
    \]
Then, we use (\ref{lb1}) with $\Omega_n=2$, and choose $R_n$ such that
    \beq
    \frac{2}{R_n}=o\Big(\frac{1}{\sqrt[2s]{k_n}}\Big),
    \eeq
so that $R_n \gg \sqrt[2s]{k_n}$. Applying (\ref{essbd}) and (\ref{rogersbound}), we see that
we need that $R_n^{k_n}=e^{o\big(\frac{n}{\sqrt[s]{k_n}}\big)}$, so that $k_n=o\big(\frac{n^{\frac{s}{s+1}}}{\log{n}}\big)$
is sufficient. Finally, by Theorem
\ref{Iuitrekenen3} and (\ref{IinftyTaylor}), we have that
    \beq
    I_{k_n}\big(1\pm \frac{1}{\sqrt[2s]{k_n}}\big)\geq \frac 1{4\sqrt[s]{k_n}}
    \Big(1+{\cal O}\big(\frac{1}{\sqrt[s]{k_n}}\big)\Big).
    \eeq
This completes the proof.
\qed

We now discuss the above results. In \cite{HHK02},
it was conjectured that when $s=2$, the rate of a {\it single} bit error
for a fixed user is asymptotic to $\frac 1{2\sqrt{k}}$ when $k\rightarrow \infty$.
See also \cite{HHK02b}. We see that we
obtain a similar result, but our constant is $1/4$ rather than the expected
$1/2$. On the other hand, our result is valid for {\it all} $s\geq 2$.

Related results where obtained for a related model, {\it Hard-Decision Parallel
Interference Cancelation} (HD-PIC) where bits are iteratively estimated by bits,
i.e., the estimates are rounded to $\pm 1$. Thus, this scheme is not linear,
as SD-PIC is. In \cite{hofstad, HK02}, similar results
as the above are obtained, and it is shown that the rate for a bit-error for a given user
is asymptotic to $\frac{s}{8}\sqrt[s]{\frac{4}{k}}$ when $s$ is fixed and
$k\rightarrow \infty$. This result is similar in spirit as the one
in Theorem \ref{thm-fsSDPIC} above. The explanation of why the rate tends to zero
as $1/\sqrt[s]{k}$ is much simpler for the case of SD-PIC, where the relation
to eigenvalues is rather direct, compared to the explanation for HD-PIC, which is much
more elaborate. It is interesting to see that both when $s=\infty$ and when
$s$ is finite and $k\rightarrow \infty$, the rates in the two systems are
of the same order.

Interestingly, in \cite{LV04}, it was shown
that for $s=1$ and $k_n=\frac{n}{\gamma\log{n}},$ with high probability, all
bits are estimated correctly when $\gamma<2$, while, with high probability, there is at
least one bit-error when $\gamma>2$. Thus, $k_n={\cal O}(\frac{n}{\log{n}})$ is
critical for the MF system, where we do not apply SD-PIC. For SD-PIC with an arbitrary
number of staged of SD-PIC, we have no
bit-errors with large probability for all $k_n=\frac{n}{\gamma \log n}$ for {\it all}
$\gamma>0$, and we can even pick larger values of $k_n$ such that $k_n=
o(\frac{n}{\log\log n})$. Thus, SD-PIC is more efficient than MF, in the sense
that it allows more users to transmit without creating bit-errors. Furthermore,
in \cite{LV04}, the results proved in
this paper are used for a further
comparison between SD-PIC, HD-PIC and MF. Unfortunately, when we only
apply a {\it finite} number of stages of SD-PIC, we can only allow for
$k_n=o\big(\frac{n^{\frac s{s+1}}}{\log{n}}\big)$ users.
Similar results were obtained for HD-PIC when $k_n=o(\frac{n}{\log{n}})$.

We close this discussion on SD-PIC and HD-PIC by noting that for $k=\beta n$,
$\lambda_{\rm min}$ converges to $(1-\sqrt\beta)_+^2$, while the
largest eigenvalue $\lambda_{\rm max}$ converges to $(1+\sqrt\beta)^2$
(see \cite{bai2, bai1, YinBaiKri88}).
This is explained in more detail in \cite{grant},
and illustrates that SD-PIC has no bit-errors with
probability converging to 1 whenever $\beta<(\sqrt{2}-1)^2\approx 0.17...$
However, unlike the case where $k_n=o(\frac{n}{\log\log n})$, we do not
obtain bounds on how the probability of a bit-error tends to zero.

A further CDMA system is the {\it decorrelator}, which explicitly
inverts the matrix $W$ (without approximating it by the partial sum
$\sum_{\varsigma=0}^{s-1} (I-W)^{\varsigma}$). One way of doing so is
to fix a large value $M$ and to compute
    \beq
    \hat{\mathbf{Z}}^{\tinbr{s}}_{\sss M}=M^{-1}\sum_{\varsigma=0}^{s-1}
    \big(I-\frac{W}{M}\big)^{\varsigma}W \mathbf{Z},
    \eeq
and
   \beq
    \hat b_{m,\sss{M}}^{\smallsup{s}} = {\rm sign}(\mathbf{\hat{Z}}^{\smallsup{s}}_{m,\sss{M}}).
    \eeq
This is a certain {\it weighted} SD-PIC scheme. This scheme will converge to
$\mathbf{b}$ as $s\rightarrow \infty$ whenever $\lambda_{\rm min}>0$ and
$\lambda_{\rm max}<M$.
By taking $M$ such that $I_{\sss \infty}(M)\geq \log{2}$, and using Proposition \ref{prop-lambdanul},
we obtain the following result:

\begin{theorem}[Bit-error rate for optimal weighted SD-PIC]
\label{thm-wSDPIC}
For all $k$ fixed, or for
$k=k_n\rightarrow \infty$ such that $k_n=
o(\frac{n}{\log\log n})$ and $M$ such that $I_{\sss \infty}(M)\geq \log{2}$,
    \beq
    -\frac {1}{n} \log \prob_k\big(\exists m=1, \ldots, k \mbox{ for which }\hat{b}_{m,\sss{M}}^{\tinbr{s}}\neq
    b_m\big) \geq \log 2.
    \eeq
\end{theorem}

\noindent
The above result can even be generalised to $k_n$ that grow arbitrarily fast with $n$,
by taking $M$ dependent on $n$. For example, when we take $M>k_n$, then
$\lambda_{\rm max}\leq k_n<M$ is guaranteed.
\vskip0.5cm

Further interesting problems arise when we allow the received signal to be
{\it noisy}. In this case, the bit-error can be caused either by the properties
of the eigenvalues, as in the case when there is no noise, or by the noise.
When there is noise, weighted SD-PIC for large $M$ {\it enhances} the noise,
which makes the problem significantly harder. See \cite{klok} for further
details. A solution to Conjecture \ref{conj-small} may prove to be useful in such
an analysis.

\paragraph{Acknowledgement.}
The work of AF and RvdH was supported in part by NWO.
The work of RvdH and MJK was performed in part at Delft University
of Technology.

\end{document}